\documentclass{amsart}
\title{Kronecker-Weber plus epsilon}
\author{Greg W. Anderson}
\address{School of Mathematics, University of Minnesota,
Minneapolis, Minnesota 55455}
\email{gwanders@math.umn.edu}
\date{This paper has been published. The journal citation is
Duke Math.\ J.\ \textbf{114}(2002), 439--475.}

\newcommand{\Pbold}{{\mathbf P}}
\newcommand{\Hbold}{{\mathbf H}}

\newcommand{\DDD}{{\mathbf D}}

\DeclareMathOperator{\Gal}{{\rm Gal}}

\DeclareMathOperator{\id}{{\rm id}}
\newcommand{\iso}{{\stackrel{\sim}{\rightarrow}}}
\newcommand{\QQ}{{\mathbb Q}}
\DeclareMathOperator{\ab}{{\rm ab}}
\newcommand{\QQbar}{\overline{\QQ}}

\newcommand{\ZZ}{{\mathbb Z}}
\newtheorem{Proposition}[subsection]{Proposition}
\newtheorem{LittleProposition}[subsubsection]{Proposition}
\newtheorem{Lemma}[subsection]{Lemma}
\newcommand{\CC}{{\mathbb C}}

\DeclareMathOperator{\sign}{{\rm sign}}

\theoremstyle{remark}
{\newtheorem{Remark}[subsubsection]{Remark}}

\newcommand{\AAA}{{\mathcal A}}
\newcommand{\ebold}{{\mathbf e}}

\newcommand{\abold}{{\mathbf a}}
\newcommand{\bbold}{{\mathbf b}}
\newcommand{\cbold}{{\mathbf c}}
\newcommand{\RR}{{\mathbb R}}

\DeclareMathOperator{\image}{{\mathrm
image}}
\newcommand{\NN}{{\mathbb N}}
\newcommand{\xbold}{{\mathbf x}}

\begin{document}
\maketitle
\begin{abstract} 
We say that a group is {\em almost abelian} if every
commutator is central and squares to the identity.
Now let $G$ be the Galois group of the algebraic closure
of the field $\QQ$ of rational numbers in the field of
complex numbers. Let
$G^{\ab+\epsilon}$ be the quotient of
$G$ universal for continuous homomorphisms to almost abelian
profinite groups and let
$\QQ^{\ab+\epsilon}/\QQ$ be the corresponding Galois
extension. We prove that
$\QQ^{\ab+\epsilon}$ is generated by the roots of unity, the
fourth roots of the rational primes and the square roots of
certain algebraic sine-monomials. The inspiration for the paper
came from recent  studies of algebraic
$\Gamma$-monomials by P.~Das and by S.~Seo. 
\end{abstract}

\section{Introduction}
        
We say that a
group is {\em almost abelian} if every
commutator is central and squares to the
identity. Let $G$ be the Galois group of the algebraic
closure of the field of rational numbers $\QQ$ in
the field $\CC$ of complex numbers. Let
$G^{\ab+\epsilon}$ be the quotient of $G$
universal for continuous homomorphisms to almost abelian
profinite groups. Let
$G^\epsilon$ be the kernel of the natural
map of $G^{\ab+\epsilon}$ to the
abelianization $G^{\ab}$ of $G$. By
construction the group
$G^{\epsilon}$ is central in
$G^{\ab+\epsilon}$ and killed by $2$.
 Let
$\QQ^{\ab}$ (resp.~$\QQ^{\ab+\epsilon}$) be the
Galois extension of $\QQ$ in $\CC$ with Galois group $G^{\ab}$
(resp.~$G^{\ab+\epsilon}$). The
Kronecker-Weber theorem determines the
structure of the group
$G^{\ab}$ and provides an explicit description of the
field
$\QQ^{\ab}$. The theory of
\cite{Frohlich}
in principle determines the structure of the group
$G^{\ab+\epsilon}$ but does not provide an 
explicit description of the field
$\QQ^{\ab+\epsilon}$. Kummer theory identifies the
Pontryagin dual of
$G^\epsilon$ with
$H^0(G^{\ab},\QQ^{\ab\times}/\QQ^{\ab\times 2})$. 
Our  purpose in this paper is to exhibit for the
latter group an explicit
$\ZZ/2\ZZ$-basis, thereby obtaining a description of
the field $\QQ^{\ab+\epsilon}$ as explicit as
that provided  for the field $\QQ^{\ab}$ by the
Kronecker-Weber theorem. Our method is 
more or less elementary and
independent of the theory of
\cite{Frohlich}.  The inspiration for our work
came from the recent studies
\cite{Das} and
\cite{Seo} of algebraic $\Gamma$-monomials.

Our main results are as follows. Let
$\AAA$ be the free abelian group on symbols of the
form  
$$[a]\;\;\;\;(a\in \QQ),$$ 
modulo the identifications
$$ [a]=[b]\Leftrightarrow a-b\in \ZZ.$$
For
all prime numbers $p<q$, if $2<p$ put
$$\abold_{pq}:=
\sum_{i=1}^{\frac{p-1}{2}}\left(\left[\frac{i}{p}\right]
-\sum_{k=0}^{\frac{q-1}{2}}\left[\frac{i}{pq}+\frac{k}{q}\right]
\right)-
\sum_{j=1}^{\frac{q-1}{2}}\left(\left[\frac{j}{q}\right]
-\sum_{\ell=0}^{\frac{p-1}{2}}\left[\frac{j}{pq}+\frac{\ell}{p}\right]
\right),$$
e.~g.,
$$\abold_{3\cdot 5}=\left[\frac{1}{3}\right]+
\left[\frac{2}{15}\right]-\left[\frac{4}{15}\right]
-\left[\frac{1}{5}\right],$$
and if $2=p$ put
$$\abold_{pq}:=
\left(\left[\frac{1}{4}\right]-\sum_{k=0}^{\frac{q-1}{2}}
\left[\frac{1}{4q}+\frac{k}{q}\right]\right)-
\sum_{j=1}^{\frac{q-1}{2}}
\left(\left[\frac{j}{q}\right]
+\left[-\frac{1}{2q}+\frac{j}{q}\right]
-\left[\frac{j}{2q}\right]
-\left[-\frac{1}{4q}+\frac{j}{2q}\right]
\right),$$
e.~g.,
$$\abold_{2\cdot 3}:=\left[\frac{1}{4}\right]
-\left[\frac{5}{12}\right]-\left[\frac{1}{3}\right].$$
Let 
$$\sin:\AAA\rightarrow \QQ^{\ab\times}$$
be the unique homomorphism such that
$$\sin [a]=\left\{\begin{array}{cl}
2\sin \pi a\left(=\left|1-e^{2\pi i
a}\right|\right)&\mbox{if
$0<a<1$}\\ 1&\mbox{if
$a=0$}
\end{array}\right.\;\;\;\;(a\in\QQ\cap[0,1)).$$
We
prove that the family of real numbers
$$
\left\{\sqrt{\ell}\right\}_{\ell:\mbox{\scriptsize
prime}}\bigcup
\left\{
\sin\abold_{pq}\right\}_{\begin{subarray}{c}
p,q:\mbox{\scriptsize prime}\\
p<q\end{subarray}}
$$
projects to a $\ZZ/2\ZZ$-basis of the group
$H^0\left(G^{\ab},\QQ^{\ab\times}/\QQ^{\ab\times
2}\right)$ and hence that 
$$\QQ^{\ab+\epsilon}=
\QQ^{\ab}\left(
\left\{\sqrt[4]{\ell}\right\}_{\ell:\mbox{\scriptsize
prime}}\bigcup
\left\{
\sqrt{\sin\abold_{pq}}\right\}_{\begin{subarray}{c}
p,q:\mbox{\scriptsize prime}\\
p<q\end{subarray}}\right).$$
We actually prove more.
We define a canonical injective homomorphism
$$\DDD:H^0\left(G^{\ab},\QQ^{\ab\times}/\QQ^{\ab\times
2}\right)\rightarrow\bigwedge^2
H^1(G^{\ab},\ZZ/2\ZZ)
\left(=\bigwedge^2 \QQ^\times/\QQ^{\times 2}\right)$$
and exhibit a preimage  for each element of a
natural $\ZZ/2\ZZ$-basis of the target. (In
particular, it turns out that $\DDD$ is an
isomorphism.) Our main results are the
``Auxiliary Formula''
(\S\ref{subsection:EasyFirstApplication}) and
the ``Main Formula''
(\S\ref{subsubsection:MainFormula}). Our main
technical tools are the ``Log Wedge Formula''
(\S\ref{subsection:LogWedge}) and a family  of
combinatorial identities of ``Das type''
(\S\ref{subsection:DasIdentities}). 

Now the fact that the homomorphism
$\DDD$ is an isomorphism is not really new: one can easily
deduce it from \cite[Theorem 4.10,
p.~56]{Frohlich}. (See
Remark~\ref{Remark:FrohlichComment} below for
further discussion.) Rather, it is the explicit
procedure developed here for inverting the map $\DDD$  that is
really new. 

A key role  is played in this paper by the
{\em universal odd
ordinary distribution}
$U^-$, namely the quotient of
$\AAA$ by the subgroup generated by all expressions of
the form
$$[a]-\sum_{i=0}^{N-1}\left[\frac{a+i}{N}\right],\;\;\;
[a]-[1-a]\;\;\;(a\in \QQ,\;\mbox{$N$: positive
integer}).$$ 
A study of
$U^-$ and related objects was made in
\cite{Kubert1} and
\cite{Kubert2}, building on \cite{Sinnott}.
In particular, it was proved that the
torsion subgroup of
$U^-$ is killed by $2$.  Kubert's
results combined with the idea behind the
algebraicity criterion of
\cite{KoblitzOgus} yield a more or less
mechanical procedure for determining whether a
given element of
$\AAA$ represents a torsion element of $U^-$.
The double complex
method of
\cite{Anderson} yields a canonical $\ZZ/2\ZZ$-basis
for the torsion subgroup of $U^-$ indexed by finite
sets of prime numbers of even cardinality. In
\cite[Sec.~3 and Sec.~9]{Das} it is proved that the family
$\{\abold_{pq}\}_{2<p<q}$ represents the ``two-odd-prime''
part of the canonical basis. The method of Das can easily be
modified to prove  that the family
$\{\abold_{pq}\}_{p<q}$ represents the ``two-prime'' part of
the canonical basis.

 The torsion subgroup $U^-$ plays an important role in the
theory of algebraic $\Gamma$-monomials.  Let
$$\Gamma:\AAA\rightarrow
\RR^\times$$ be the unique homomorphism such that
$$\Gamma \left([a]\right)=
\left\{\begin{array}{cl}
\sqrt{2\pi}/\Gamma(a)&\mbox{if
$0<a<1$}\\ 1&\mbox{if $a=0$}
\end{array}\right.\;\;\;\;\;\;\;(a\in \QQ\cap [0,1)).
$$
Now fix
$\abold\in
\AAA$ representing a torsion element of $U^-$.
By  straightforward manipulation of standard
functional equations satisfied by $\Gamma(s)$ one
verifies that
$\Gamma(\abold)$ is an algebraic number. 
Numbers of the form $\Gamma(\abold)$ are the so called
{\em algebraic $\Gamma$-monomials}.
 Under mild hypotheses, a reciprocity law
\cite[Thm.~7.15, p.~91]{DMOS} due to Deligne 
links $\Gamma(\abold)$ to an analogously
defined Jacobi sum Hecke character. Deligne reciprocity
is the main motivation for studying
algebraic $\Gamma$-monomials.

In \cite{Das} the algebraic
$\Gamma$-monomials were related to
algebraic sine-monomials and properties of the
latter were investigated.  Again
fix
$\abold\in
\AAA$ representing a torsion element of $U^-$.
By
\cite[Thm.~6]{Das}, under mild hypotheses,
$\Gamma(\abold)$ factors  as the root of a
rational number   times
$\sqrt{\sin\abold}$, and the factorization can in
principle be worked out explicitly. For example, one has
$$
\Gamma(\abold_{3\cdot 5})
=\frac{\Gamma\left(\frac{4}{15}\right)
\Gamma\left(\frac{1}{5}\right)}{\Gamma\left(\frac{1}{3}\right)
\Gamma\left(\frac{2}{15}\right)}
=3^{-\frac{1}{5}}5^{\frac{1}{12}}\sqrt{\frac{\sin\frac{\pi}{3}
\cdot\sin\frac{2\pi}{15}}
{\sin\frac{4\pi}{15}\cdot
\sin\frac{\pi}{5}}}=3^{-\frac{1}{5}}5^{\frac{1}{12}}
\sqrt{\sin\abold_{3\cdot 5}}
$$
and
$$
\Gamma(\abold_{2\cdot 3})=\frac{
\Gamma\left(\frac{5}{12}\right)\Gamma\left(\frac{1}{3}\right)}
{\sqrt{2\pi}\,\Gamma\left(\frac{1}{4}\right)}
=2^{-\frac{1}{4}}3^{\frac{1}{8}}
\sqrt{\frac{\sin\frac{\pi}{4}}{2\cdot\sin
\frac{5\pi}{12}\cdot
\sin
\frac{\pi}{3}}}
=2^{-\frac{1}{4}}3^{\frac{1}{8}}\sqrt{\sin
\abold_{2\cdot 3}}.
$$
Indeed, for
all primes
$p<q$ one has
$$\begin{array}{rcl}
\displaystyle\frac{\Gamma(\abold_{pq})}{\sqrt{\sin\abold_{pq}}}
&=&\displaystyle
\left\{\begin{array}{cl}
p^{\frac{-(q-1)^2}{16q}}q^{\frac{(p-1)^2}{16p}}
&\mbox{if $2<p$,}\\\\
2^{-\frac{q-1}{8}}q^{\frac{1}{8}}&\mbox{if $2=p$,}
\end{array}\right.
\end{array}$$
as is explained below in 
Remark~\ref{Remark:GammaMonomials}.
By
\cite[Thm.~11]{Das}, under mild hypotheses, the extension
$\QQ^{\ab}\left(\sqrt{\sin\abold}\right)/\QQ$ is
Galois, and in fact the latter is true in general; see
Proposition~\ref{Proposition:GaloisProperty} below. 
 By \cite[Thm.~22]{Das}, if
$\abold$ represents an element of the canonical basis
for the torsion subgroup of $U^-$ indexed by a finite
set of odd primes of cardinality at least four, then
$\sqrt{\sin \abold}\in \QQ^{\ab}$.

The theory of \cite{Das} begs the question
as to the nature of the Galois extensions of the
form $\QQ^{\ab}\left(\sqrt{\sin\abold_{pq}}\right)/\QQ$.
In order to begin answering that question, Das
calculated (see \cite[end of Sec.~16]{Das})
 the structure of the
central  extension
$$1\rightarrow \{\pm 1\}
\rightarrow \Gal\left(\QQ\left(e^{2\pi
i/60},\sqrt{\sin\abold_{3\cdot 5}}\right)/\QQ\right)
\rightarrow (\ZZ/60\ZZ)^\times \rightarrow 1$$
explicitly and in particular proved that the
middle group is nonabelian.   That calculation of Das
was the primary inspiration for this paper.

Another important inspiration was the
striking result
\cite[Prop.~2.3]{Seo} according to which, for
all odd primes $p<q$, one has
$$(-1)^{v_q(\sin\abold_{pq})}=\left(\begin{array}{c}
p\\
q\end{array}\right),\;\;\;
(-1)^{v_p(\sin\abold_{pq})}=\left(\begin{array}{c}
q\\
p\end{array}\right)$$
where $\left(\begin{array}{c}
\cdot\\
\cdot\end{array}\right)$ is the Legendre symbol,
$v_p$ is any additive valuation of
$\QQ^{\ab}$ above $p$ such that
$v_p\left(1-e^{\frac{2\pi i}{p}}\right)=1$, and $v_q$
is analogously chosen above $q$.  It follows that if
$\left(\begin{array}{c} p\\
q\end{array}\right)=-1$ or $\left(\begin{array}{c}
q\\
p\end{array}\right)=-1$, then 
$\sqrt{\sin\abold_{pq}}$ cannot belong to $\QQ^{\ab}$.
Seo's result greatly encouraged the author
to attempt the calculations detailed in this paper.

\section{Group-theoretical background}

\subsection{The category of $\NN$-profinite groups}
We say that a topological group $G$ is {\em $\NN$-profinite}
if $G$ is compact and
there exists a chain
$U_1\supseteq U_2\supseteq U_3\supseteq \dots$
of open normal subgroups of $G$ forming a neighborhood base
at the origin. A {\em morphism}
of
$\NN$-profinite groups is by definition a continuous homomorphism
of topological groups; thus the $\NN$-profinite groups
form a category. The Galois group of
any Galois extension of countable fields is
$\NN$-profinite.  (But the Galois group of an
algebraic closure of the field of rational functions in one
variable over the field of complex numbers is {\em not}
$\NN$-profinite.) 
A compact topological group $G$
is $\NN$-profinite if and only if $G$ has a countable
neighborhood base consisting of open compact sets.

\begin{Proposition}\label{Proposition:ProfiniteSplitting}
Let $G$ be an $\NN$-profinite group. Let $H$ be a closed
normal subgroup of $G$. There exists a function
$\vartheta:G\rightarrow G$ with the following properties:
\begin{itemize}
\item $\vartheta$ is continuous.
\item $\vartheta$ is
constant on cosets of
$H$ in
$G$.  
\item $\vartheta(\sigma)\in \sigma H$ for all $\sigma\in G$.
\end{itemize}
\end{Proposition}
\proof Fix a descending chain $U_1\supseteq U_2\supseteq
U_3\supseteq \dots$ of open normal subgroups of $G$ forming a
neighborhood base at the origin. 
Choose a set $S_1$ of representatives for
the cosets of $HU_1$ in $G$. For each index
$i>1$ choose a set
$S_i\subseteq U_{i-1}$ of representatives for the cosets of
$HU_i$ in
$HU_{i-1}$. 
For each index $i$ let $\vartheta_i:G\rightarrow G$
be the unique map collapsing each coset of $HU_i$ in $G$ to its
unique representative in the finite set $S_1\cdots S_i$.
The uniform limit
$\vartheta$ of the maps
$\vartheta_i$ exists and has all the desired properties.
\qed

\subsection{Cochains, cocycles, coboundaries and cohomology
classes} Let $G$ be a topological group.
Let $A$ be a topological abelian group equipped with a
continuous left $G$-action. 
  For
each nonnegative integer $n$,  let
${C}^n(G,A)$ denote the group of continuous functions
$$a=((\sigma_1,\dots,\sigma_n)\mapsto
a_{\sigma_1,\dots,\sigma_n}):G^n\rightarrow A,$$
and equip the graded
group  
$${C}^*(G,A):=\bigoplus{C}^n(G,A)$$ with a differential
$\delta$ of degree $1$ by the standard rule
$$(\delta a)_{\sigma_1,\dots,\sigma_{n+1}}=
\sigma_1a_{\sigma_2,\dots,\sigma_{n+1}}
+\cdots+(-1)^i
a_{\sigma_1,\dots,\sigma_i\sigma_{i+1},\dots,\sigma_{n+1}}
+\cdots+(-1)^{n+1}a_{\sigma_1,\dots,\sigma_n}.$$
Elements of $C^n(G,A)$ will be called {\em cochains} of
degree
$n$. 
In low degree 
 the differential
$\delta$ takes the form
$$\begin{array}{rcll}
(\delta a)_{\sigma_1}&=&\sigma_1 a-a&\mbox{if $n=0$,}\\
(\delta a)_{\sigma_1,\sigma_2}
&=&\sigma_1
a_{\sigma_2}-a_{\sigma_1\sigma_2}+a_{\sigma_1},&\mbox{if
$n=1$,}\\
(\delta a)_{\sigma_1,\sigma_2,\sigma_3}
&=&\sigma_1
a_{\sigma_2,\sigma_3}-a_{\sigma_1\sigma_2,\sigma_3}
+a_{\sigma_1,\sigma_2\sigma_3}-a_{\sigma_1,\sigma_2}&\mbox{if
$n=2$.}
\end{array}
$$
Elements of the groups
$$\begin{array}{rcl}
Z^n(G,A)&:=&\ker\left(\delta:C^n(G,A)\rightarrow
C^{n+1}(G,A)\right),\\
B^n(G,A)&:=&\image\left(\delta:C^{n-1}(G,A)\rightarrow
C^n(G,A)
\right),\\
H^n(G,A)&:=&Z^n(G,A)/B^n(G,A)
\end{array}$$ will be called
{\em $n$-cocycles}, {\em $n$-coboundaries}
and {\em cohomology classes of degree $n$},
respectively. If $G$ is profinite and $A$ is discretely
topologized, then
$H^n(G,A)$ is the usual Galois cohomology
group of degree $n$.

\subsection{Cup product}
Let $G$ be a topological group and let $A$ be a commutative
topological ring equipped with a continuous $G$-action. 
 Given
$a\in C^p(G,A)$ and $b\in C^q(G,A)$,
put
$$(a\cup b)_{\sigma_1,\dots,\sigma_{p+q}}=
a_{\sigma_1,\dots,\sigma_p}\cdot\sigma_1\cdots
\sigma_p b_{\sigma_{p+1},\dots,\sigma_{p+q}},$$
thereby defining the {\em cup product} 
$a\cup b\in C^{p+q}(G,A)$
of the cochains $a$ and $b$.
For $p=q=1$ one has
$$(a\cup
b)_{\sigma_1,\sigma_2}=a_{\sigma_1}\cdot
\sigma_1 b_{\sigma_2}.$$
In general one has
$$\delta(a\cup b)=(\delta a)\cup b+(-1)^p a\cup(\delta
b)$$
and thus the cup product construction induces a product
in $H^*(G,A)$.

\subsection{Extensions in the category of $\NN$-profinite
groups} 

Let
$G$ be an $\NN$-profinite group and let $A$ be an
$\NN$-profinite abelian group equipped with a continuous left
$G$-action. An exact sequence
$$\Sigma:\;\;\;1\rightarrow A\stackrel{i}{\rightarrow}U
\stackrel{p}{\rightarrow} G\rightarrow 1$$
in the category of $\NN$-profinite groups such that
$$ui(a)u^{-1}=i(p(u)a)\;\;\;(u\in U,\;a\in A)$$
is called an {\em extension}
of $G$ by $A$ in the category of $\NN$-profinite groups. 
If $A$ is a trivial
$G$-module, the homomorphism $i:A\rightarrow U$ takes values in
the center of $U$; in such a case one says that
$\Sigma$ is a {\em central} extension. Extensions
$\Sigma$ and
$\Sigma'$ of $G$ by $A$ are said to be {\em isomorphic} if there
exists an isomorphism $\psi:U\iso U'$ of $\NN$-profinite
groups such that
$\psi\circ i=i'$ and $p'\circ \psi =p$.
A {\em set-theoretic splitting}
$\vartheta:G\rightarrow U$ of $\Sigma$ is a continuous map
such that $p\circ
\vartheta=\id_G$.
Set-theoretic splittings exist by
Proposition~\ref{Proposition:ProfiniteSplitting}. We say that
$\Sigma$ is  {\em trivial} if there exists a
set-theoretic splitting that is also a group homomorphism.
Any two trivial extensions of $G$ by $A$ are isomorphic. Each
set-theoretic splitting
$\vartheta$ of
$\Sigma$ gives rise to a cochain $a\in C^2(G,A)$  by the rule
$$\vartheta(\sigma)\vartheta(\tau)=
i(a_{\sigma,\tau})\vartheta(\sigma\tau) \;\;\;(\sigma,\tau\in
G).$$
Associativity of the group law in $U$
forces $a$ to be a $2$-cocycle.
The cohomology class of the
$2$-cocycle
$a$ depends only the isomorphism class of the extension
$\Sigma$ and accordingly is denoted $[\Sigma]$. Every
$2$-cocycle representing the cohomology class $[\Sigma]$ arises
from some set-theoretic splitting of $\Sigma$. The construction
$\Sigma\mapsto [\Sigma]$ puts isomorphism classes of extensions
of $G$ by $A$ in the category of $\NN$-profinite groups into
bijective correspondence with $H^2(G,A)$ and sends every
trivial extension to $0$.

\subsection{Almost abelian $\NN$-profinite groups}
Let $G$ be an $\NN$-profinite group and 
put
$$[\sigma,\tau]:=\sigma\tau\sigma^{-1}\tau^{-1}\;\;\;
(\sigma,\tau\in G).$$
The {\em abelianization} $G^{\ab}$ of $G$ is the abelian
quotient of $G$ universal for morphisms to abelian
$\NN$-profinite groups. One has 
$$G^{\ab}=G/[G,G],$$
where $[G,G]$ is the closed subgroup of $G$
topologically generated by the set $$\{[\sigma,\tau]\mid
\sigma,\tau\in G\}.$$ We say that $G$
is {\em almost abelian} if $[\sigma,\tau]$ is central and
$2$-torsion for all $\sigma,\tau\in G$.
 We define $G^{\ab+\epsilon}$ to be the almost abelian quotient of
$G$ universal for morphisms to almost abelian $\NN$-profinite
groups. One has
$$G^{\ab+\epsilon}=G/[G,G]^\epsilon,$$
where $[G,G]^\epsilon$ is the closed subgroup of $[G,G]$
topologically generated by the set
$$\{[\sigma,\tau]^2\mid \sigma,\tau\in
G\}\cup\{[\sigma,[\tau,\eta]]\mid \sigma,\tau,\eta\in G\}.$$
Put
$$G^\epsilon:=[G,G]/[G,G]^\epsilon.$$
The preceding constructions fit together to form a canonical
central extension
$$\Sigma_G^\epsilon:\;\;1\rightarrow G^\epsilon
\rightarrow G^{\ab+\epsilon}
\rightarrow G^{\ab}\rightarrow 1$$
of $\NN$-profinite groups.  We remark that the Pontryagin dual of
$G^\epsilon$ is canonically isomorphic to
$H^1(G^\epsilon,\ZZ/2\ZZ)$  and that the latter is a vector
space over $\ZZ/2\ZZ$ of at most countably infinite dimension.

\begin{Proposition}\label{Proposition:DDefinition3}
Let $G$ be an $\NN$-profinite group.
 Let
$H^2(G^{\ab},\ZZ/2\ZZ)^-$ be the quotient of
$H^2(G^{\ab},\ZZ/2\ZZ)$
by the subgroup consisting of all cohomology classes
of the form $[\Sigma]$ for some extension $\Sigma$ of
$G^{\ab}$ by $\ZZ/2\ZZ$ the middle group of which is
abelian. For each $c\in H^1(G^\epsilon,\ZZ/2\ZZ)$,
view $c$ as a continuous homomorphism $G^\epsilon\rightarrow 
\ZZ/2\ZZ$, and let
$$(a\mapsto c\circ a):H^2(G^{\ab},G^\epsilon)\rightarrow
H^2(G^{\ab},\ZZ/2\ZZ)$$ be the map  induced by $c$ in
degree
$2$ cohomology. The homomorphism
$$(c\mapsto
c\circ[\Sigma_G^\epsilon]):H^1(G^\epsilon,\ZZ/2\ZZ)\rightarrow
H^2(G^{\ab},\ZZ/2\ZZ)$$ 
followed by the projection
$$H^2(G^{\ab},\ZZ/2\ZZ)\rightarrow H^2(G^{\ab},\ZZ/2\ZZ)^-$$
is injective.
\end{Proposition}
\proof 
Fix
$0\neq c\in H^1(G^\epsilon,\ZZ/2\ZZ)$ and consider the
extension
$$\Sigma:\;0\rightarrow \ZZ/2\ZZ\rightarrow G^{\ab+\epsilon}/\ker
c\rightarrow G^{\ab}\rightarrow 1.$$
Then
$$[\Sigma]=c\circ
[\Sigma_G^\epsilon]$$
and the middle group of $\Sigma$ is
nonabelian by definition of $G^{\ab}$.
 \qed

\begin{Lemma}\label{Lemma:Cocycling}
Let $G=\prod_{i=1}^\infty G_i$
be an abelian
$\NN$-profinite group decomposed as a countable product of
cyclic profinite groups $G_i$. Let an extension
$$\Sigma:\;\;0\rightarrow
\ZZ/2\ZZ\rightarrow U
\rightarrow G\rightarrow 1$$
of $\NN$-profinite groups be given. Let
$$a\in Z^2(G,\ZZ/2\ZZ)$$ be any $2$-cocycle representing
the  cohomology class
$$[\Sigma]\in H^2(G,\ZZ/2\ZZ).$$
For each index $i$ fix
$\sigma_i\in G$ with the following properties:
\begin{itemize}
\item $\sigma_i$ projects to a topological generator of
$G_i$.
\item $\sigma_i$ projects to $1\in G_j$ for all indices $j\neq i$.
\end{itemize}
 Put
$$\alpha_{ij}:\equiv
a_{\sigma_i,\sigma_j}+a_{\sigma_j,\sigma_i}\bmod{2}$$
for all
indices $i<j$. Then the following  hold:
\begin{enumerate}
\item $\alpha_{ij}$
vanishes for all but finitely many pairs of indices.

\item $\alpha_{ij}$ vanishes
if and only if $\sigma_i$
and $\sigma_j$ have liftings to $U$ that commute.

\item $\alpha_{ij}$ depends only on the
cohomology class $[\Sigma]$.

\item If $\sigma_i$ or $\sigma_j$ is the square of an element
of $G$, then
$\alpha_{ij}$ vanishes.
\item The middle group $U$ of $\Sigma$ is abelian if and only if
all the $\alpha_{ij}$ vanish.
\end{enumerate}
\end{Lemma}
\proof  
Fix a set-theoretic splitting $\vartheta$ of $\Sigma$
giving rise to the $2$-cocycle $a$.
 
1. By hypothesis one has
$$(\delta a)_{1,1,\sigma}-(\delta
a)_{\sigma,1,1}\equiv
a_{1,\sigma}+a_{\sigma,1}\equiv
0\bmod{2}\;\;\;(\sigma\in G)$$ and
there exists an open
normal subgroup
$V\subseteq G$ such that
$a_{\sigma,\tau}$ depends only on the pair of cosets
$(\sigma V,\tau V)$. 
It follows that if either $\sigma_i\in V$
or $\sigma_j\in V$, then $\alpha_{ij}$ vanishes. But one
has
$\sigma_i\in V$ for all but finitely many indices $i$.

2. The following statements are equivalent:
\begin{itemize}
\item The coefficient $\alpha_{ij}$ vanishes.
\item $\vartheta(\sigma_i)$ and $\vartheta(\sigma_j)$
commute.
\item $\sigma_i$ and $\sigma_j$ have liftings to $U$
that commute.
\item Any liftings of $\sigma_i$ and $\sigma_j$ to $U$
commute.
\end{itemize}

3. By what we have already proved, the coefficients
$\alpha_{ij}$ depend only on the isomorphism class of the
extension $\Sigma$ and hence only on the cohomology class
$[\Sigma]$.

4. 
Suppose, say, that 
$\sigma_i=\tau^2$ for some
$\tau\in G$. Then
$$\vartheta(\tau)^2c=\vartheta(\sigma_i),\;\;\;
\vartheta(\tau)\vartheta(\sigma_j)=\vartheta(\sigma_j)\vartheta(\tau)
c'$$ where $c$ and $c'$ belong to the center of $U$,
and hence 
$$\vartheta(\sigma_i)\vartheta(\sigma_j)=
\vartheta(\tau)\vartheta(\sigma_j)\vartheta(\tau)c'c
=\vartheta(\sigma_j)\vartheta(\tau)c'\vartheta(\tau)c'c
=\vartheta(\sigma_j)\vartheta(\sigma_i),
$$
i.~e., $\vartheta(\sigma_i)$ and
$\vartheta(\sigma_j)$ commute, and hence the coefficient
$\alpha_{ij}$ vanishes.

5. All the coefficients
$\alpha_{ij}$ vanish if and only if 
the elements of the family $\{\vartheta(\sigma_i)\}$
commute among themselves. But the union of the latter family
with $\ker(U\rightarrow G)$
generates $U$ topologically.
\qed

\begin{Proposition}\label{Proposition:DDefinition1}
Let $G$ be an abelian $\NN$-profinite group admitting 
decomposition as a countably infinite product of cyclic
profinite groups. The cup product induces an isomorphism
$$\bigwedge^2 H^1(G,\ZZ/2\ZZ)\iso
H^2(G,\ZZ/2\ZZ)^-$$ of vector spaces over $\ZZ/2\ZZ$,
where $H^2(G,\ZZ/2\ZZ)^-$ is as defined in
Proposition~\textup{\ref{Proposition:DDefinition3}}.
\end{Proposition}
\proof 
For each $e\in Z^1(G,\ZZ/2\ZZ)=H^1(G,\ZZ/2\ZZ)$, the
cohomology class of the 
$2$-cocycle $e\cup e$ vanishes in the quotient
$H^2(G,\ZZ/2\ZZ)^-$ by Lemma~\ref{Lemma:Cocycling}. 
Therefore the cup product induces a well defined map
$$\bigwedge^2 H^1(G,\ZZ/2\ZZ)\rightarrow
H^2(G,\ZZ/2\ZZ)^-.$$
Fix a decomposition $G=\prod_{i=1}^\infty G_i$ of $G$ as a
countably infinite product of cyclic profinite groups
$G_i$ and for each index
$i$ fix
$\sigma_i\in G$ with the properties specified in 
Lemma~\ref{Lemma:Cocycling}. 
Let $I$ be the set of indices $i$ such that $G_i\not\subseteq
G^2$. For each index 
$i\in I$ there exists a unique $1$-cocycle  
$$e_i\in
Z^1(G,\ZZ/2\ZZ)=H^1(G,\ZZ/2\ZZ)$$ such that
$$(e_i)_{\sigma_j}\equiv \delta_{ij}\bmod{2}$$ for all indices
$j$. Then the family
$$\{e_i\}_{i\in I}$$ is a $\ZZ/2\ZZ$-basis for
$H^1(G,\ZZ/2\ZZ)$ and hence the family 
$$\{e_i\wedge e_j\}_{\begin{subarray}{c}
i,j\in I\\
i<j\end{subarray}}$$ is a
$\ZZ/2\ZZ$-basis for
$\bigwedge^2H^1(G,\ZZ/2\ZZ)$.  By Lemma~\ref{Lemma:Cocycling}
the cohomology classes of the $2$-cocycles belonging to the
family
$$\{e_i\cup e_j\}_{\begin{subarray}{c}
i,j\in I\\
i<j\end{subarray}}$$ remain independent over $\ZZ/2\ZZ$ in
$H^2(G,\ZZ/2\ZZ)^-$.
Therefore the map in question is injective. For any
$a\in Z^2(G,\ZZ/2\ZZ)$ the cohomology class of the
two-cocycle 
$$a+\sum_{\begin{subarray}{c}
i,j\in I\\
i<j\end{subarray}}(a_{\sigma_i,\sigma_j}+a_{\sigma_j,\sigma_i})
e_i\cup e_j$$
vanishes in the quotient $H^2(G,\ZZ/2\ZZ)^-$ by
Lemma~\ref{Lemma:Cocycling}. Therefore the map in question is
surjective.
\qed

\begin{Proposition}\label{Proposition:DDefinition2}
For any $\NN$-profinite group $G$, the group
$[G,G]^\epsilon$ is the
intersection of all relatively open subgroups of $[G,G]$
that are of index $2$ in $[G,G]$ and  normal in $G$.
\end{Proposition}
\proof Let $H$ be the
intersection of all relatively open subgroups of $[G,G]$
that are of index $2$ in $[G,G]$ and  normal in $G$. 
Since the quotient
$G/H$ is almost abelian, one has $[G,G]^\epsilon\subseteq H$.
Since the continuous
linear functionals $G^\epsilon\rightarrow
\ZZ/2\ZZ$ separate points, one has
$H\subseteq [G,G]^\epsilon$.
\qed

\pagebreak

\section{The homomorphism $\DDD$}
\label{section:DDefinition}
\subsection{The setting for the rest of the paper}

\subsubsection{}
Let
$\QQbar$ be the algebraic closure of the field
$\QQ$ of rational numbers in the field $\CC$ of complex numbers.
Let $\QQ^{\ab}$ be the compositum of all subfields
$K\subseteq\QQbar$ abelian over $\QQ$. 
Let $\QQ^{\ab+\epsilon}$ be the compositum of
all subfields $K\subseteq \QQbar$ that are quadratic over
$\QQ^{\ab}$ and Galois over $\QQ$. (In the terminology of
\cite{Das} such extensions $K/\QQ$ would be called  {\em double
coverings} of $\QQ^{\ab}$.)
Let
$G$ denote the Galois group of $\QQbar/\QQ$. One has
$$G^{\ab}=\Gal(\QQ^{\ab}/\QQ),\;\;\;G^{\ab+\epsilon}=\Gal(\QQ^{\ab+\epsilon}/\QQ),\;\;\;
G^\epsilon=\Gal(\QQ^{\ab+\epsilon}/\QQ^{\ab})$$
by
Proposition~\ref{Proposition:DDefinition2} and the definitions.

\subsubsection{}
Put
$$S:=\{-1\}\cup \{\mbox{rational prime numbers}\}.$$
The family
$$\{p\bmod{\QQ^{\times 2}}\}_{p\in S}$$
is a $\ZZ/2\ZZ$-basis for $\QQ^\times/\QQ^{\times 2}$.
For each $p\in S$ we define a $1$-cocycle
$$e_p\in
Z^1(G^{\ab},\ZZ/2\ZZ)=H^1(G^{\ab},\ZZ/2\ZZ)$$
by the rule 
$$\sigma \sqrt{p}=(-1)^{(e_p)_\sigma}\sqrt{p}\;\;\;(\sigma\in
G^{\ab}).$$ 
By Kummer theory the family 
$$\{e_p\}_{p\in S}$$
is a $\ZZ/2\ZZ$-basis for
$H^1(G^{\ab},\ZZ/2\ZZ)$.

\subsubsection{}
Let $G_{-1}\subset G^{\ab}$ be the subgroup
generated  by the restriction of complex conjugation to
$\QQ^{\ab}$. For each odd prime
$p$, let $G_p\subset G^{\ab}$ be the inertia subgroup at $p$.  
Let
$G_2\subset G^{\ab}$ be the subgroup of the inertia subgroup at
$2$ fixing
$\sqrt{-1}$. 
For all $p\in S$ the profinite group $G_p$ is cyclic and has
nonzero
$2$-rank. By the Kronecker-Weber theorem the family 
$\{G_p\}_{p\in S}$ is the family of ``coordinate axes''
for a decomposition of $G^{\ab}$ into a product of
cyclic profinite groups indexed by
$S$.  It follows by
Proposition~\ref{Proposition:DDefinition1} that we can
identify
$H^2(G^{\ab},\ZZ/2\ZZ)^-$ with $\bigwedge^2
H^1(G^{\ab},\ZZ/2\ZZ)$.  
Note that the family
$$\{e_p\wedge e_q\}_{\begin{subarray}{c}
p,q\in S\\
p<q
\end{subarray}}$$
is a $\ZZ/2\ZZ$-basis for$\bigwedge^2
H^1(G^{\ab},\ZZ/2\ZZ)$.

\subsubsection{}
For each $p\in S$ we  arbitrarily fix a topological
generator
$\sigma_p$ of $G_p$.  
Note that $\sigma_{-1}$
is the restriction to $\QQ^{\ab}$ of complex 
conjugation. Note that
$$(e_p)_{\sigma_q}\equiv
\delta_{pq}\bmod{2}\;\;\;(p,q\in S),$$
i.~e., the families $\{\sigma_p\}_{p\in S}$
and $\{e_p\}_{p\in S}$ are in a convenient way dual.

\subsection{The definition of $\DDD$}

\label{subsection:DDefinition}
 
We define
$$\DDD:H^0\left(G^{\ab},\QQ^{\ab\times}/\QQ^{{\ab}\times
2}\right)\rightarrow
\bigwedge^2H^1(G^{\ab},\ZZ/2\ZZ)$$
to be the injective homomorphism obtained
as the composite of
the following three homomorphisms:
\begin{itemize}
\item The canonical isomorphism 
$$H^0(G^{\ab},\QQ^{\ab\times}/\QQ^{\ab
\times 2})\iso H^1(G^\epsilon,\ZZ/2\ZZ)
$$
provided by Kummer theory combined with
Proposition~\ref{Proposition:DDefinition2}.\\
\item The canonical injective homomorphism
$$H^1(G^\epsilon,\ZZ/2\ZZ)\rightarrow
H^2(G^{\ab},\ZZ/2\ZZ)^-$$ provided by
Proposition~\ref{Proposition:DDefinition3}.\\
\item The inverse of the canonical isomorphism
$$\bigwedge^2
H^1(G^{\ab},\ZZ/2\ZZ)\iso H^2(G^{\ab},\ZZ/2\ZZ)^-$$
provided by Proposition~\ref{Proposition:DDefinition1}.\\
\end{itemize}
The map $\DDD$ is the focus of our investigation.

\subsection{The Log Wedge Formula}
\label{subsection:LogWedge}
We carry out the
diagram-chase necessary to make the map $\DDD$ usably
explicit.

\subsubsection{}
Fix
$u\in \QQ^{\ab\times}$ such that
$$(u\bmod{\QQ^{\ab\times
2}})\in H^0(G^{\ab},\QQ^{\ab\times}/\QQ^{\ab\times 2}).
$$ 
The
$1$-cocycle 
$$c\in
Z^1(G^\epsilon,\ZZ/2\ZZ)=H^1(G^\epsilon,\ZZ/2\ZZ)$$ defined by
the rule
$$\sigma \sqrt{u}=(-1)^{c_\sigma}\sqrt{u}\;\;\;(\sigma\in
G^\epsilon)$$ is the one
corresponding via Kummer theory  to the class
$u\bmod{\QQ^{\ab\times 2}}$. 

\subsubsection{}
Let 
$$\vartheta:G^{\ab}\rightarrow G^{\ab+\epsilon}$$
be any set-theoretic splitting of the
canonical central extension $\Sigma^\epsilon_G$.
The $2$-cocycle 
$$a\in Z^2(G^{\ab},\ZZ/2\ZZ)$$
defined  by the formula
$$a_{\sigma,\tau}=c_{\vartheta(\sigma)\vartheta(\tau)
\vartheta(\sigma\tau)^{-1}}\;\;\;(\sigma,\tau\in G^{\ab})$$
represents the cohomology class $c\circ[\Sigma^\epsilon_G]\in
H^2(G^{\ab},\ZZ/2\ZZ)$. 

\subsubsection{}
Let
$$\{v_\sigma\}_{\sigma\in G^{\ab}}$$ be a locally constant
family of elements of
$\QQ^{\ab \times}$ such that
$$\sigma u= v_\sigma^2 u \;\;\;(\sigma\in G^{\ab}).$$
 Define a $1$-cochain
$$f\in C^1(G^{\ab},\ZZ/2\ZZ)$$ by the rule
$$\vartheta(\sigma) \sqrt{u}=(-1)^{f_\sigma}v_\sigma
\sqrt{u}\;\;\;(\sigma\in G^{\ab}).$$ A straightforward if
tedious calculation gives the formula
$$(-1)^{a_{\sigma,\tau}}\sqrt{u}=\vartheta(\sigma)\vartheta(\tau)\vartheta(\sigma\tau)
^{-1}\sqrt{u}=
(-1)^{(\delta f)_{\sigma,\tau}}(\delta
v)_{\sigma,\tau}\sqrt{u}\;\;\;(\sigma,\tau\in G^{\ab}).$$
The latter says that the
$2$-cocycle 
$$\delta\log_{-1} v\in
Z^2(G^{\ab},\ZZ/2\ZZ)$$ also
represents the cohomology class $c\circ[\Sigma_G^\epsilon]$,
where
$$\log_{-1}:\{\pm 1\}\iso \ZZ/2\ZZ$$
is the evident (and only) isomorphism.

\subsubsection{}
Put
$$\alpha_{p,q}:=\frac{(\delta v)_{\sigma_p,\sigma_q}}{
(\delta v)_{\sigma_q,\sigma_p}}=\frac
{\sigma_pv_{\sigma_q}/v_{\sigma_q}}
{\sigma_q v_{\sigma_p}/v_{\sigma_p}}\in \{\pm 1\}\;\;\;
\left(\begin{array}{c}
p,q\in S\\
p<q\end{array}\right).$$  By
Lemma~\ref{Lemma:Cocycling} 
one has $\alpha_{p,q}=1$
for all but finitely many pairs $\{p<q\}$ and the image
in
$H^2(G^{\ab},\ZZ/2\ZZ)^-$ of the cohomology
class of the $2$-cocycle $$
\delta \log_{-1}v+\sum_{\begin{subarray}{c}
p,q\in S\\
p<q
\end{subarray}}\log_{-1}\alpha_{p,q}\cdot
e_p\cup e_q$$ 
vanishes. 
We arrive finally at the result
$$\DDD(u\bmod{\QQ^{\ab\times 2}})=
\sum_{\begin{subarray}{c}
p,q\in S\\
p<q\end{subarray}}\log_{-1}\alpha_{p,q}\cdot e_p\wedge
e_q$$
to which in the sequel we refer as the {\em Log Wedge
Formula}. 

\subsection{The Auxiliary Formula}
\label{subsection:EasyFirstApplication}
Here is an easy first application of the Log Wedge Formula. Let
$\ell$ be any prime number. One has
$\sqrt{\ell}\in\QQ(\zeta_{4\ell})^\times$ and hence
there exists a unique function
$$(\sigma\mapsto v_\sigma):G^{\ab}\rightarrow \{1,i\}$$
factoring through $\Gal(\QQ(\zeta_{4\ell})/\QQ)$
such that
$$\sigma
\sqrt{\ell}=v_\sigma^2\sqrt{\ell}\;\;\;\left(\mbox{equivalently:
}v_\sigma^{2}=(-1)^{(e_\ell)_\sigma}\right)\;\;\;
\;\;\;\;(\sigma\in G^{\ab}).$$
For all $p,q\in S$ one has
$$\sigma_pv_{\sigma_q}/v_{\sigma_q}=\left\{\begin{array}{rl}
-1&\mbox{if $p=-1 $ and $q=\ell$,}\\
1&\mbox{otherwise,}
\end{array}\right.
$$
and hence
$$\DDD \left(\sqrt{\ell}\bmod{\QQ^{\ab\times 2}}\right)
=e_{-1}\wedge e_\ell
$$
by the Log Wedge
Formula. We refer to
the result above as the {\em Auxiliary Formula}.

\section{The Main Formula}
  
\subsection{The universal ordinary distribution
and related apparatus}
\subsubsection{} We denote by $\AAA$
the free abelian group on symbols of the form
$$[a]\;\;\;\;(a\in \QQ)$$
modulo the identifications
$$[a]=[b]\Leftrightarrow a-b\in \ZZ.$$
Put
$$\ebold(a):=e^{2\pi i a}\;\;\;\;(a\in \QQ).$$
We equip
$\AAA$ with the unique action of
$G^{\ab}$  such that
$$\sigma[a]=[b]\Leftrightarrow \sigma
\ebold(a)=\ebold(b)\;\;\;(a,b\in \QQ
,\;\sigma\in
G^{\ab}).$$  

\subsubsection{}
For each prime
number
$p$ we define
$$Y_p:\AAA\rightarrow \AAA$$
to be the unique endomorphism such that
$$Y_p[a]:=[a]-\sum_{i=0}^{p-1}
\left[\frac{a+i}{p}\right]\;\;\;(a\in \QQ).$$
The operator $Y_p$ is $G^{\ab}$-equivariant. 
The operators $Y_p$ commute among themselves.
Put
$$U:=\frac{\AAA}{\sum_p Y_p \AAA},\;\;\;U^-:=
\frac{\AAA}{(1+\sigma_{-1})\AAA+\sum_p Y_p\AAA},$$
thereby defining the {\em universal ordinary distribution}
and the {\em universal odd ordinary distribution},
respectively. The action of  $G^{\ab}$ on $\AAA$
descends to 
$U$ and
$U^-$. By \cite{Kubert1}, the group $U$
is free abelian. It follows that the torsion subgroup of
$U^-$ can be identified with the cohomology group
$H^2(\langle \sigma_{-1}\rangle,U)$, and in particular
that the torsion subgroup of $U^-$ is killed by $2$. By
\cite{Kubert2}, the group $G^{\ab}$ operates trivially on
the torsion subgroup of $U^-$. 

\subsubsection{}
We denote by $\AAA'$ the
subgroup of $\AAA$ generated by all symbols of the form 
$$[a]\;\;\;\;\left(a\in \QQ\setminus
\frac{1}{2}\ZZ\right).$$
 Notice that $\AAA'$ is stable
under the action of $G^{\ab}$ and stable also under the action
of all the operators $Y_p$.
We call a homomorphism
$$H:\AAA'\rightarrow \AAA'$$
a {\em lifting operator} if the identities
$$H^2=H,\;\;\;H+\sigma_{-1}H\sigma_{-1}=1$$
hold, in which case
$$\AAA'=H\AAA'\oplus\sigma_{-1}H\AAA'$$ and one has implications
$$\begin{array}{rcl}
(1+\sigma_{-1})\xbold=0&\Rightarrow&
\xbold=(1-\sigma_{-1})H\xbold\\
 (1-\sigma_{-1})
\xbold=0&\Rightarrow&
\xbold=(1+\sigma_{-1})H\xbold
\end{array}\;\;\;\;\;\;(\xbold\in \AAA').
$$

\subsubsection{}
To each partition
$$(0,1)\cap
\left(\QQ\setminus
\frac{1}{2}\ZZ\right)=T\coprod \{1-a\mid a\in T\}$$ is
attached a unique lifting operator
$$\Hbold_T:\AAA'\rightarrow \AAA'$$
such that
$$\Hbold_T[a]=\left\{\begin{array}{cl}
[a]&\mbox{if $a\in T$}\\
0&\mbox{if $a\not\in T$}
\end{array}\right.\;\;\;\;\left(a\in (0,1)\cap
\left(\QQ\setminus \frac{1}{2}\ZZ\right)\right).$$ Following
\cite{Das}, we call
$$\Hbold:=\Hbold_{(0,1/2)\cap \QQ}$$
the {\em canonical lifting operator}. 

\subsubsection{}
We define 
$$\xi:\AAA\rightarrow \QQ^{\ab\times}$$
to be the unique homomorphism such that
$$\xi [a]=\left\{\begin{array}{cl}
1-\ebold(a)&\mbox{if $a\neq 0$}\\
1&\mbox{if $a=0$}
\end{array}\right.\;\;\;(a\in [0,1)\cap \QQ).$$
The homomorphism $\xi$ is Galois
equivariant. 
One
has
$$\xi\left(Y_p\left[a\right]\right)=\left\{\begin{array}{cl}
1&\mbox{if $a\neq 0$}\\ p^{-1}=(p^{-1/2})^2&\mbox{if $a=
0$}
\end{array}\right.\;\;\;\;(\mbox{$p$: prime}, 
\;a\in [0,1)\cap \QQ)$$ 
and
$$\xi\left((1+\sigma_{-1})[a]\right)=\left\{\begin{array}{cl}
1&\mbox{if $a=0$}\\
\left|1-\ebold(a)\right|^2=4\sin^2\pi a&\mbox{if
$a\neq 0$}
\end{array}\right.
\;\;\;\;(a\in [0,1)\cap \QQ).$$
It follows that $\xi$ induces a  $G^{\ab}$-equivariant
homomorphism 
$$U^-\rightarrow
\QQ^{\ab\times}/\QQ^{\ab\times 2}$$
and hence a homomorphism
$$(\mbox{torsion
subgroup of $U^-$})\rightarrow
H^0\left(G^{\ab},\QQ^{\ab\times}/\QQ^{\ab\times 2}\right).$$
Note also that $\xi$ kills $\sum_p Y_p\AAA'$.

\subsubsection{}\label{subsubsection:DeviousParityOperator}
We define a homomorphism
$$\sin:\AAA\rightarrow \QQ^{\ab\times}$$
by the rule
$$\sin [a]:=\left\{\begin{array}{cl}
2\sin \pi a&\mbox{if $a\in \QQ\cap (0,1)$,}\\
1&\mbox{if $a=0$.}
\end{array}\right.$$
One then has
$$\sin \abold=|\xi(\abold)|\;\;\;(\abold\in \AAA).$$
Let $\AAA''$ denote the subgroup of $\AAA'$ generated by elements
of the form 
$$[a]\;\;\;\;(a\in (0,1)\cap \QQ\cap \ZZ_2),$$
where $\ZZ_2$ denotes the $2$-adic completion of $\ZZ$.  Let 
$$\Pbold:\AAA''\rightarrow \AAA''$$ be the unique
idempotent endomorphism such that 
$$\Pbold[a]=\left\{\begin{array}{cl}
[a]&\mbox{if $a$ is a $2$-adic unit}\\
0&\mbox{otherwise}
\end{array}\right.\;\;\;\;(a\in (0,1)\cap
\QQ\cap\ZZ_2).$$
Now $a\in (0,1)\cap \QQ\cap \ZZ_2$
 is a $2$-adic unit if and only if $1-a$ is not a $2$-adic
unit; it follows that one has an identity
$$\Pbold\abold+\sigma_{-1}\Pbold\sigma_{-1}\abold=\abold\;\;\;(\abold\in
\AAA'').$$
Let 
$$\deg:\AAA\rightarrow \ZZ$$ be the
unique homomorphism such that 
$$\deg [a]=1\;\;\;\;(a\in \QQ).$$
\begin{LittleProposition}\label{Proposition:SignRule}
One has
$$\sign \sigma\sin \abold= \frac{\sigma
\sin \abold}{\sin \sigma\abold}= (-1)^{\deg
\Pbold(1-\sigma)\abold}\left(
\frac{\sigma\sqrt{-1}}
{\sqrt{-1}}\right)^{\deg
\abold}\;\;\; (\abold\in
\AAA'',\;\sigma\in G^{\ab}).$$ 
\end{LittleProposition}
\proof For each $a\in \QQ\cap\ZZ_2$ 
there exists a solution $\tilde{a}\in \QQ\cap \ZZ_2$
of the congruence 
$$2\tilde{a}\equiv
a\bmod{\ZZ}$$
unique modulo $\ZZ$ because $(\QQ\cap \ZZ_2)/\ZZ$
is the subgroup of $\QQ/\ZZ$ consisting of elements of odd
order. For all $a\in \QQ\cap (0,1)\cap \ZZ_2$ one has
$$\frac{a}{2}\equiv \tilde{a}+\frac{1}{2}\bmod{\ZZ}
\Leftrightarrow \mbox{$a$ is a $2$-adic unit.}
$$
It follows that one has an identity
$$|1-\ebold(a)|=2\sin \pi a=(-1)^{\deg
\Pbold[a]}\;\;\frac{\ebold(\tilde{a})-
\ebold(-\tilde{a})}{\ebold(1/4)}\;\;\;(a\in
(0,1)\cap \QQ\cap \ZZ_2).$$ Let
$$\tilde{\xi}:\AAA''\rightarrow
\QQ^{\ab\times}$$ be the unique homomorphism such that
$$\tilde{\xi}[a]=\ebold(\tilde{a})-\ebold(-\tilde{a})\;\;\;
(a\in
(0,1)\cap \QQ\cap \ZZ_2).$$
Then
$\tilde{\xi}$ is
$G^{\ab}$-equivariant and one has an identity
$$\sin \abold=(-1)^{\deg
\Pbold\abold}\frac{\tilde{\xi}(\abold)}
{\ebold(1/4)^{\deg
\abold}}\;\;\;
(\abold\in \AAA''),$$ whence the result.
\qed

\begin{LittleProposition}\label{Proposition:GaloisProperty}
For every
$\abold\in
\AAA$ representing a torsion element of $U^-$, the extension
$\QQ^{\ab}(\sqrt{\sin\abold})/\QQ$ is Galois. 
(Cf.~\cite[Thm.~11]{Das}.)
\end{LittleProposition}
\proof For every $\abold\in \AAA$
the ratio $\xi(\abold)/\sin\abold$
is a root of unity. Therefore the homomorphism
$$(\mbox{torsion subgroup of $U^-$})\rightarrow
H^0\left(G^{\ab},\QQ^{\ab\times}/\QQ^{\ab\times 2}\right)$$
induced by $\xi:\AAA\rightarrow
\QQ^{\ab\times}$ is also induced by
$\sin:\AAA\rightarrow \QQ^{\ab\times}$.
\qed

\subsection{Identities of Das type}
\label{subsection:DasIdentities}
We carry out by  {\em ad hoc} methods tailored to our needs
what are in effect diagram-chases through the double complex
employed by
\cite{Das}.

\subsubsection{}
\label{subsubsection:FirstDasIdentity}
Fix a lifting
operator
$$H:\AAA'\rightarrow
\AAA',$$ 
 distinct prime numbers $p$ and $q$, and $\xbold\in \AAA$ such
that
$$Y_p\xbold,Y_q\xbold\in \AAA',\;\;\;\;(1-\sigma_{-1})\xbold=0.$$
One has
$$(1-\sigma_{-1})Y_p\xbold=0,\;\;\;(1-\sigma_{-1})Y_q\xbold=0$$
and hence
$$Y_p\xbold=(1+\sigma_{-1})HY_p\xbold,\;\;\;Y_q\xbold=(1+\sigma_{-1})HY_q\xbold.$$
One has
$$(1+\sigma_{-1})(Y_qHY_p-Y_pHY_q)\xbold=
(Y_q(1+\sigma_{-1})HY_p-Y_p(1+\sigma_{-1})HY_q)\xbold=0,$$
hence
$$(1-\sigma_{-1})(HY_pHY_q-HY_qHY_p)\xbold=
(Y_pHY_q-Y_qHY_p)\xbold,$$
and hence
$$\begin{array}{cl}
&2(HY_pHY_q-HY_qHY_p)\xbold\\\\
=&(Y_pHY_q-Y_qHY_p)\xbold
+(1+\sigma_{-1})(HY_pHY_q-HY_qHY_p)\xbold.
\end{array}$$
We refer to the latter relation as the {\em first Das
identity}.

\subsubsection{}
Fix lifting operators
$$H,\bar{H}:\AAA'\rightarrow \AAA',$$
distinct prime numbers
$p$ and
$q$, and $\xbold\in \AAA$ such that
$$Y_p\xbold,Y_q\xbold\in \AAA',\;\;\;\;(1-\sigma_{-1})\xbold=0.$$
Put
$$\abold:=(HY_pHY_q-HY_qHY_p)\xbold,\;\;\;
\bar{\abold}:=(\bar{H}Y_p\bar{H}Y_q-\bar{H}Y_q\bar{H}Y_p)\xbold,$$
and
$$\bbold:=(Y_p H(H-\bar{H})Y_q-Y_qH(H-\bar{H})Y_p)\xbold,\;\;
\cbold:=H(\abold-\bar{\abold})-H\bbold.$$
One has
$$(1+\sigma_{-1})(H-\bar{H})Y_q\xbold=0,\;\;\;
(1+\sigma_{-1})(H-\bar{H})Y_p\xbold=0,$$
hence
$$\begin{array}{rcl}
(H-\bar{H})Y_q\xbold&=&(1-\sigma_{-1})H(H-\bar{H})Y_q\xbold,\\
(H-\bar{H})Y_p\xbold&=&(1-\sigma_{-1})H(H-\bar{H})Y_p\xbold,
\end{array}$$
 hence
$$(1-\sigma_{-1})(\abold-\bar{\abold})=
(Y_p(H-\bar{H})Y_q-Y_q(H-\bar{H})Y_p)\xbold=(1-\sigma_{-1})\bbold$$
and hence
$$\abold-\bar{\abold}=\bbold+(1+\sigma_{-1})\cbold.
$$
We refer to the latter relation as the {\em second
Das identity}.

\subsubsection{}
Fix a lifting operator 
$$H:\AAA'\rightarrow \AAA',$$
 distinct
primes $p$ and $q$, and $\xbold\in \AAA$ such that
$$Y_p\xbold,Y_q\xbold\in
\AAA',\;\;\;(1-\sigma)\xbold=0\;\;\;(\sigma\in G^{\ab}).$$ 
Put
$$\abold:=(HY_pHY_q-HY_qHY_p)\xbold.$$
For each $\sigma\in G^{\ab}$ put\\
$$\begin{array}{rcll}
\bbold_\sigma&:=&
Y_pH(1-\sigma)HY_q\xbold-Y_qH(1-\sigma)HY_p\xbold&\in
Y_p\AAA'+Y_q\AAA'\\\\
\cbold_\sigma
&:=& H(1-\sigma)\abold-H\bbold_\sigma&\in \AAA'.
\end{array}$$
After making the substitution
$\bar{H}=\sigma
H\sigma^{-1}$
in the second Das identity
and simplifying, one obtains the family
of identities
$$(1-\sigma)\abold=\bbold_\sigma+(1+\sigma_{-1})\cbold_\sigma
\;\;\;(\sigma\in G^{\ab})$$
to which we refer in the sequel as the {\em Das conjugation
formula}. Note
that all the expressions
$\bbold_\sigma$ belong to the kernel of the homomorphism
$\xi$. It follows that the Das conjugation formula
specializes under the homomorphism $\xi$ to a family of
numerical identities
$$\sigma\xi(\abold)=\xi(\abold)/\sin^{2}\cbold_\sigma
\;\;\; (\sigma\in G^{\ab})$$
in which the expressions $\bbold_\sigma$ figure not at all.
\begin{Remark}
We got the idea for the technical tool that
we are calling the Das conjugation formula
from the intriguing example \cite[end of Sec.~16]{Das}. 
\end{Remark}

\subsection{Das classes and the statement of the Main
Formula}

\subsubsection{} Let a lifting operator
$$H:\AAA'\rightarrow \AAA'$$ and primes $p<q$ be given.
Put
$$\abold:=(HY_pHY_q-HY_qHY_p)\cdot\left\{\begin{array}{cl}
\;[0]&\mbox{if $2<p<q$,}\\
\;[0]+[1/2]&\mbox{otherwise.}
\end{array}\right.$$
In the present situation the
identities of Das type have the following implications:
\begin{itemize}
\item $\abold$ represents a
$2$-torsion $G^{\ab}$-invariant element of 
$U^-$ independent of $H$.
\item 
$\xi(\abold)$ is a real number and hence
$\xi(\abold)=\pm \sin \abold$.
\item $\sin \abold$ represents a class in
$H^0\left(G^{\ab},\QQ^{\ab\times}/\QQ^{\ab\times 2}\right)$
independent of 
$H$. 
\item One has a family of numerical identities
$$\sigma \sin
\abold =\sin
\abold/\sin^2\cbold_\sigma\;\;\; (\sigma\in G^{\ab})$$ 
where  $\{\cbold_\sigma\}$ is 
the cochain figuring in the Das conjugation formula.
\end{itemize}
We call the torsion
element of
$U^-$ represented by
$\abold$ the {\em Das class} associated to the pair $\{p<q\}$
of prime numbers.  
\subsubsection{}
For all primes $p<q$, put
$$\abold_{pq}:=(\Hbold Y_p\Hbold Y_q-\Hbold Y_q \Hbold
Y_p)\cdot\left\{\begin{array}{cl}
\;[0]&\mbox{if $2<p<q$,}\\
\;[0]+[1/2]&\mbox{otherwise,}
\end{array}\right.$$
thereby defining the {\em canonical
representative} of the Das class associated to the
pair $\{p<q\}$ of primes. For $2<p<q$ one has
$$\abold_{pq}=
\left(\sum_{i=1}^{\frac{p-1}{2}}\left(\left[\frac{i}{p}\right]
-\sum_{k=0}^{\frac{q-1}{2}}\left[\frac{i}{pq}+\frac{k}{q}\right]
\right)\right)-
\left(\sum_{j=1}^{\frac{q-1}{2}}\left(\left[\frac{j}{q}\right]
-\sum_{\ell=0}^{\frac{p-1}{2}}\left[\frac{j}{pq}+\frac{\ell}{p}\right]
\right)\right).$$ 
For $2=p<q$ one has
$$\abold_{pq}=\left(
\left[\frac{1}{4}\right]-
\sum_{k=0}^{\frac{q-1}{2}}
\left[\frac{1}{4q}+\frac{k}{q}\right]\right)-
\sum_{j=1}^{\frac{q-1}{2}}
\left(\left[\frac{j}{q}\right]
+\left[-\frac{1}{2q}+\frac{j}{q}\right]
-\left[\frac{j}{2q}\right]
-\left[-\frac{1}{4q}+\frac{j}{2q}\right]
\right).$$
The discovery that $\abold_{pq}$ for $2<p<q$
represents  a torsion element of $U^-$ is due to
Das; see \cite[Sec.~9]{Das}.

\begin{Remark}\label{Remark:GammaMonomials}
The
double complex method of
\cite{Anderson} gives rise to a canonical $\ZZ/2\ZZ$-basis
for the torsion subgroup of $U^-$ indexed by finite sets
of prime numbers of even cardinality.
In \cite[\S3 and \S9]{Das} it is proved that the family
$\{\abold_{pq}\}_{2<p<q}$ represents the ``two-odd-prime''
part of the canonical basis.
The method of Das can easily be modified to show that
the family
$\{\abold_{pq}\}_{p<q}$ represents ``two-prime'' part 
of the canonical basis. 
\end{Remark}

\subsubsection{}
\label{subsubsection:MainFormula}The main result
of this paper  is the
relation
$$\DDD\left(\sin \abold_{pq}\bmod{\QQ^{\ab\times
2}}\right) =e_p\wedge e_q\;\;\;(p<q\mbox{ 
prime})$$ to which we refer in the sequel as the {\em
Main Formula}.  The Main and Auxiliary Formulas
together immediately imply that
$$\QQ^{\ab+\epsilon}=
\QQ^{\ab}\left(
\left\{\sqrt[4]{\ell}\right\}_{\ell:\mbox{\scriptsize
prime}}\bigcup
\left\{
\sqrt{\sin\abold_{pq}}\right\}_{\begin{subarray}{c}
p,q:\mbox{\scriptsize prime}\\
p<q\end{subarray}}\right).$$
Thus we get a quite explicit description of the
field $\QQ^{\ab+\epsilon}$.

\begin{Remark}\label{Remark:FrohlichComment}
The following is a corollary to the Main and
Auxiliary Formulas:
\begin{quote}
$(\star)$ For $G=\Gal(\QQbar/\QQ)$, the canonical injective map
$$(c\mapsto 
c\circ[\Sigma_G^\epsilon]):H^1(G^\epsilon,\ZZ/2\ZZ)\rightarrow
H^2(G^{\ab},\ZZ/2\ZZ)^-$$ 
defined in Proposition~\textup{\ref{Proposition:DDefinition3}} is
an isomorphism.
\end{quote}
But $(\star)$ is not really new. One can get it very easily
from the theory of \cite{Frohlich}, as we now explain. Let
$S$ be a finite set of primes of
$\QQ$ including 
$2$ and the infinite prime. 
 Let $\Gamma$ be the maximal pro-$2$ quotient of $G$ unramified
outside the primes in $S$.
To prove $(\star)$ it is enough to prove
that the canonical injective map
$$(e\mapsto
e\circ[\Sigma_{\Gamma}^\epsilon]):H^1(\Gamma^\epsilon,\ZZ/2\ZZ)
\rightarrow H^2(\Gamma^{\ab},\ZZ/2\ZZ)^-$$
is an isomorphism. 
Let  $(\tau\mapsto \tau^{\ab}):\Gamma\rightarrow \Gamma^{\ab}$
be the canonical projection. For each odd finite prime $p\in S$,
choose any $\tau_p\in \Gamma$ such that $\tau^{\ab}_p$
topologically generates the inertia subgroup of $\Gamma^{\ab}$ at $p$.
Choose any $\tau_2\in \Gamma$
such that $\tau_2^{\ab}$  topologically generates the 
subgroup of the inertia 
subgroup of $\Gamma^{\ab}$ at $2$ fixing $\sqrt{-1}$. 
Choose any $\tau_{\infty}\in \Gamma$ so that
$\tau^{\ab}_\infty\in \Gamma^{\ab}$ is the automorphism
induced by complex conjugation.
The family $\{\tau_p\}_{p\in S}$ is then a minimal set of 
generators for $\Gamma$; in particular,
the $2$-rank of $\Gamma^{\ab}$ is $N$, where $N$ is the cardinality
of $S$.
By \cite[Theorem 4.10, p.~56]{Frohlich} a complete set
of defining relations for $\Gamma$ as a pro-$2$ group
is lifted from the family of relations
$$\left(\tau_p^{\ab}\right)^{p-1}=1\;\;\;(p:\mbox{odd finite prime in $S$});
\;\;\; \left(\tau^{\ab}_\infty\right)^2=1.$$
It follows that the family of
commutators
$$\{[\tau_p,\tau_q]\}_{\begin{subarray}{c}
p,q\in S\\
p<q
\end{subarray}}
$$ 
projects to a $\ZZ/2\ZZ$-basis for $\Gamma^{\epsilon}$;
in particular, the
$2$-rank of $\Gamma^{\epsilon}$ is $\frac{N(N-1)}{2}$.
By Proposition~\ref{Proposition:DDefinition1},
the $2$-rank of $H^2(\Gamma^{\ab},\ZZ/2\ZZ)^-$
is $\frac{N(N-1)}{2}$. By dimension-counting it follows that the
injective map
$e\mapsto e\circ[\Sigma_{\Gamma}^\epsilon]$ is bijective. Thus
$(\star)$ is proved {\em
\`{a} la} Fr\"{o}hlich.
\end{Remark}

\begin{Remark}
Let
$\Gamma:\AAA\rightarrow
\RR^\times$ be the unique homomorphism such that
$$\Gamma\left( [a]\right)=
\left\{\begin{array}{cl}
\sqrt{2\pi}/\Gamma(a)&\mbox{if
$0<a<1$}\\ 1&\mbox{if $a=0$}
\end{array}\right.\;\;\;\;\;\;\;(a\in \QQ\cap [0,1)).
$$
The well known functional equations satisfied
by the
$\Gamma$-function become in the present context the
identities
$$\Gamma([a]+[-a])=2\sin
\pi a,\;\;\;\;\Gamma(Y_p[a])=p^{1/2-a}\;\;\;(a\in
\QQ\cap(0,1)).$$ 
For all primes $p<q$, we then have
$$\begin{array}{rcl}
\displaystyle\frac{\Gamma(\abold_{pq})}{\sqrt{\sin\abold_{pq}}}
&=&\displaystyle\sqrt{\Gamma\left((Y_p\Hbold Y_q-Y_q\Hbold
Y_p)\cdot
\left\{\begin{array}{cl}
\,[0]&\mbox{if $2<p$}\\
\,[0]+[1/2]&\mbox{if $2=p$}
\end{array}\right)\right.}\\\\
&=&\displaystyle
\left\{\begin{array}{cl}
p^{\frac{-(q-1)^2}{16q}}q^{\frac{(p-1)^2}{16p}}
&\mbox{if $2<p$}\\\\
2^{-\frac{q-1}{8}}q^{\frac{1}{8}}&\mbox{if $2=p$}
\end{array}\right.
\end{array}$$
by the first Das identity.
The formula above (in the case
$2<p$) was discovered by Das (see
\cite[Sec.~9]{Das}) and the proof of that formula is where we
got the idea for the technical tool that we are
calling the first Das identity.
\end{Remark}

\begin{Remark}
The relations standing between the Main Formula,
the index formulas of \cite{Sinnott},
Deligne reciprocity
\cite[Thm.~7.15, p.~91]{DMOS}, the
theory of
\cite{Frohlich}, the theory of
\cite{Das},  the theory of the group cohomology of the
universal ordinary distribution (see \cite{Ouyang} and
references therein) and 
Stark's conjecture and its variants (see \cite{Tate}) deserve
to be thoroughly investigated. We have only scratched the
surface here. Stark's conjecture is relevant in view of the
well known expansion
$$\sum_{n=0}^\infty
\frac{1}{(n+x)^s}
=\frac{1}{2}-x+s\log \frac{\Gamma(x)}{\sqrt{2\pi}}+O(s^2)$$
of the Hurwitz zeta
function at $s=0$.
\end{Remark}

\begin{Remark}
The papers \cite{Thakur},
\cite{Sinha}, and \cite{Bae} (just to mention the
first several coming to mind)
suggest that there ought to be an analogue of the  Main
Formula over a global field of positive characteristic
equipped with a distinguished place.
\end{Remark}
\begin{Remark}
Perhaps there is an analogue of the Main Formula over an imaginary
quadratic field involving elliptic units. This possibility seems
especially intriguing.
\end{Remark}

\section{Proof of the Main Formula}
Fix prime numbers $p<q$.
\subsection{Reductions}
\subsubsection{}
Put
$$\xbold:=\left\{\begin{array}{cl}
\;[0]&\mbox{if $2<p$,}\\
\;[0]+[1/2]&\mbox{if $2=p$.}
\end{array}\right.
$$ 
Fix a partition
$$(0,1)\cap\left(\QQ\setminus\frac{1}{2}\ZZ\right)=
T\coprod \{1-a\mid a\in T\}$$
and put
$$H:=\Hbold_T.$$
Presently we are going to make an advantageous choice for the
set $T$. Put
$$\abold:=\left(H Y_pH Y_q-
H Y_qH Y_p\right)\xbold.$$
Then $\abold$ represents the Das class associated to
the pair $\{p<q\}$ of primes
and
$$\sin\abold_{pq}\equiv\sin\abold\bmod{\QQ^{\ab\times
2}}.$$

\subsubsection{} In the present situation the Das conjugation
formula reads
$$(1-\sigma)\abold=
\bbold_\sigma+(1+\sigma_{-1})\cbold_\sigma
\;\;\;(\sigma\in G^{\ab})$$
where
$$
\begin{array}{rcl}
\bbold_\sigma&:=&
\left(Y_pH(1-\sigma)H Y_q -Y_q
H(1-\sigma)H Y_p\right)\xbold,\\\\
\cbold_\sigma
&:=&H(1-\sigma)\abold-H\bbold_\sigma,
\end{array}$$
and one has a family of numerical identities
$$\sigma \sin \abold
=\sin \abold/\sin^2\cbold_\sigma\;\;\;(\sigma\in
G^{\ab}).$$
By construction the
function
$$\left(\sigma\mapsto\sin\cbold_\sigma\right):G^{\ab}\rightarrow
\QQ^{\ab\times }$$ 
factors 
through
$\Gal\left(\QQ\left(\ebold\left(\frac{1}{2pq}\right)\right)/
\QQ\right)$ and takes values in
$\QQ\left(\ebold\left(\frac{1}{4pq}\right)\right)$.

\subsubsection{}
The Log Wedge Formula says that 
$$\DDD\left(\sin\abold_{pq}\bmod{\QQ^{\ab\times
2}}\right)=\DDD\left(\sin\abold\bmod{\QQ^{\ab\times
2}}\right)=\sum_{\begin{subarray}{l} r,s\in S\\
r<s\end{subarray}}
\log_{-1}\alpha_{rs}\cdot e_r\wedge e_s$$
where
$$\alpha_{rs}:=
\frac{\sigma_s
\sin\cbold_{\sigma_r}/\sin\cbold_{\sigma_r}}
{\sigma_r
\sin\cbold_{\sigma_s}/\sin \cbold_{\sigma_s}}\in \{\pm 1\}.$$
(We employ $r$ and $s$ as dummy variables here since $p$ and
$q$ have already been fixed.) Since all numbers of the form
$\sin
\abold$ are positive real, we can calculate
$\alpha_{rs}$ just by keeping track of signs, i.~e., we have
$$\alpha_{rs}=\frac{\sign\sigma_s
\sin\cbold_{\sigma_r}}
{\sign\sigma_r
\sin\cbold_{\sigma_s}}.
$$
One has $\sin\cbold_{\sigma_{-1}}=1$ and hence
$\alpha_{-1,\ell}=1$ for all primes
$\ell$. For all
primes
$\ell$ distinct from $p$ and $q$, one has
$\sin\cbold_{\sigma_\ell}=\sin\cbold_1=1$ 
and moreover $\sigma_\ell$ operates trivially on the field
$\QQ\left(\ebold\left(\frac{1}{4pq}\right)\right)$.
(Recall
that $\sigma_2 \sqrt{-1}=\sqrt{-1}$.) It follows that
$\alpha_{rs}=1$ for all primes $r<s$ such that
$\{r,s\}\neq
\{p,q\}$.  
It remains now only to prove that $\alpha_{pq}=-1$.

\subsection{Calculation of $\alpha_{pq}$ 
in the case $2<p$} 
Assume  that $2<p$, i.~e., that both $p$ and $q$ are odd.
\subsubsection{}
We have
$$\alpha_{pq}=(-1)^{\deg\left( 
\Pbold(1-\sigma_q)\cbold_{\sigma_p}-\Pbold(1-\sigma_p)
\cbold_{\sigma_q}\right)}$$
by Proposition~\ref{Proposition:SignRule}. 
In this way we reduce the calculation of $\alpha_{pq}$
to an essentially combinatorial problem.

\subsubsection{}
Put
$$\Theta_p:=\sigma_p^0+\cdots+\sigma_p^{\frac{p-3}{2}},\;\;\;\;\;\;
\Theta_q:=\sigma_q^0+\cdots+\sigma_q^{\frac{q-3}{2}},$$
thereby defining elements of the integral group ring of
$G^{\ab}$. The relations
$$\left.\begin{array}{l}
(1-\sigma_p)\Theta_p[a+b]=
\left(1-\sigma_p^{\frac{p-1}{2}}\right)[a+b]=
[a+b]-[-a+b]\\\\
(1-\sigma_q)\Theta_q[a+b]=
\left(1-\sigma_q^{\frac{q-1}{2}}\right)[a+b]=
[a+b]-[a-b]
\end{array}\right\}\;\;\;\;\;
\left(a\in \frac{1}{p}\ZZ,\;\;b\in \frac{1}{q}\ZZ\right)$$
figure crucially in our calculations.

\subsubsection{}
We make an advantageous choice for the set
$T$ now. We  assume that
$T$ has been chosen in such a way that the associated lifting
operator
$H=\Hbold_T$ has the following properties:
\\
$$\begin{array}{rcl}
H\sigma^i_p\left[\frac{q}{p}\right]&=&
\left\{\begin{array}{cl}
\sigma^i_p\left[\frac{q}{p}\right]&\mbox{for $0\leq
i<(p-1)/2$,}\\\\0 &\mbox{for $(p-1)/2\leq
i<p-1$,}
\end{array}\right.\\\\
H\sigma^i_q\left[\frac{p}{q}\right]&=&
\left\{\begin{array}{cl}
\sigma^i_q\left[\frac{p}{q}\right]&\mbox{for $0\leq
i<(q-1)/2$,}\\\\
0&\mbox{for $(q-1)/2\leq
i<q-1,$}
\end{array}\right.\\\\
H\sigma_p^i\sigma_q^j\left[-\frac{1}{p}+\frac{1}{q}\right]
&=&\left\{\begin{array}{cl}
\sigma_p^i\sigma_q^j\left[-\frac{1}{p}+\frac{1}{q}\right]
&
\mbox{for $0\leq i<(p-1)/2$ and $0\leq j<q-1$,}\\\\
0&\mbox{for
$(p-1)/2\leq i<p-1$ and $0\leq j<q-1$.}
\end{array}\right.
\end{array}
$$
The reason for making this
choice is to force a lot of cancellation. 
It was in order to have the freedom to make this choice
that we took the trouble to
formulate the definition of the Das class in terms of an
arbitrarily chosen lifting operator. 

\subsubsection{} 
We make some preliminary calculations.
One has
$$Y_q\left[0\right]=-\Theta_q
\left(\left[\frac{p}{q}\right]+
\left[-\frac{p}{q}\right]\right),\;\;\;
Y_p\left[0\right]=-\Theta_p
\left(\left[\frac{q}{p}\right]+
\left[-\frac{q}{p}\right]\right),$$
hence
$$HY_q\left[0\right]=
-\Theta_q\left[\frac{p}{q}\right],\;\;\;
HY_p\left[0\right]= -\Theta_p
\left[\frac{q}{p}\right],$$
and hence
$$H(1-\sigma_q)H
Y_q\left[0\right]=-\left[\frac{p}{q}\right],\;\;\;
H(1-\sigma_p)HY_p\left[0\right]=-\left[\frac{q}{p}\right].
$$
One then has
$$\begin{array}{rcl}
Y_pH(1-\sigma_q)H
Y_q\left[0\right]&=&-
\left(\left[\frac{p}{q}\right]-\left[\frac{1}{q}\right]
-\Theta_p\left(
\left[\frac{1}{p}+\frac{1}{q}\right]+
\left[-\frac{1}{p}+\frac{1}{q}\right]\right)\right),\\\\
Y_qH(1-\sigma_p)H
Y_p\left[0\right]&=&-
\left(\left[\frac{q}{p}\right]-\left[\frac{1}{p}\right]
-\Theta_q\left(
\left[\frac{1}{p}+\frac{1}{q}\right]
+\left[\frac{1}{p}-\frac{1}{q}\right]\right)\right),\\\\
Y_pH
Y_q\left[0\right]&=&-\Theta_q
\left(\left[\frac{p}{q}\right]-\left[\frac{1}{q}\right]
-\Theta_p\left(
\left[\frac{1}{p}+\frac{1}{q}\right]+
\left[-\frac{1}{p}+\frac{1}{q}\right]\right)\right),\\\\
Y_qH
Y_p\left[0\right]&=&-\Theta_p
\left(\left[\frac{q}{p}\right]-\left[\frac{1}{p}\right]
-\Theta_q\left(
\left[\frac{1}{p}+\frac{1}{q}\right]
+\left[\frac{1}{p}-\frac{1}{q}\right]\right)\right).
\end{array}
$$
\subsubsection{}
It follows  that\\
$$
\begin{array}{rcl}
\abold&=&(HY_pHY_q-HY_qHY_p)[0]\\\\&=&
\Theta_p
\Theta_q\left[-\frac{1}{p}+\frac{1}{q}\right]
-H\Theta_q
\left(\left[\frac{p}{q}\right]-\left[\frac{1}{q}\right]
\right)+H\Theta_p
\left(\left[\frac{q}{p}\right]-\left[\frac{1}{p}\right]
\right),\\\\
H(1-\sigma_p)\abold&
=&\Theta_q
\left[-\frac{1}{p}+\frac{1}{q}\right]
+H(1-\sigma_p)H\Theta_p
\left(\left[\frac{q}{p}\right]-\left[\frac{1}{p}\right]
\right),\\\\
HY_qH(1-\sigma_p)H
Y_p\left[0\right]&=&
-H
\left(\left[\frac{q}{p}\right]-\left[\frac{1}{p}\right]
\right),\\\\

\cbold_{\sigma_p}&=&H(1-\sigma_p)\abold-
HY_pH(1-\sigma_p)HY_q[0]+HY_qH(1-\sigma_p)HY_p[0]\\\\
&=&\Theta_q
\left[-\frac{1}{p}+\frac{1}{q}\right]
+H(1-\sigma_p)H\Theta_p
\left(\left[\frac{q}{p}\right]-\left[\frac{1}{p}\right]
\right)-
H\left(\left[\frac{q}{p}\right]-
\left[\frac{1}{p}\right]\right),\\\\
H(1-\sigma_q)\abold
&=&\Theta_p
\left(\left[-\frac{1}{p}+\frac{1}{q}\right]-
\left[-\frac{1}{p}-\frac{1}{q}\right]\right)
-H(1-\sigma_q)H\Theta_q
\left(\left[\frac{p}{q}\right]-
\left[\frac{1}{q}\right]\right),
\\\\
HY_pH(1-\sigma_q)H
Y_q\left[0\right]&=&
\Theta_p\left[-\frac{1}{p}+\frac{1}{q}\right]-H
\left(\left[\frac{p}{q}\right]-\left[\frac{1}{q}\right]\right)
,\\\\
\cbold_{\sigma_q}&
=&H(1-\sigma_q)\abold-HY_pH(1-\sigma_q)HY_q[0]+
HY_qH(1-\sigma_q)HY_p[0]\\\\
&=&-\Theta_p
\left[-\frac{1}{p}-\frac{1}{q}\right]
-H(1-\sigma_q)H\Theta_q
\left(\left[\frac{p}{q}\right]-\left[\frac{1}{q}\right]
\right)+H\left(\left[\frac{p}{q}\right]-
\left[\frac{1}{q}\right]\right).
\end{array}$$

\subsubsection{}
Finally, one has
$$ 
\begin{array}{cl}
&\displaystyle
\Pbold(1-\sigma_q)\cbold_{\sigma_p}-\Pbold(1-\sigma_p)
\cbold_{\sigma_q}\\\\
=&\displaystyle\Pbold\left(
\left[-\frac{1}{p}+\frac{1}{q}\right]
-\left[-\frac{1}{p}-\frac{1}{q}\right]
\right)-\Pbold\left(-\left[-\frac{1}{p}-\frac{1}{q}\right]
+\left[\frac{1}{p}-\frac{1}{q}\right]\right)\\\\
=&\displaystyle \Pbold
\left[-\frac{1}{p}+\frac{1}{q}\right]
-\Pbold\left[\frac{1}{p}-
\frac{1}{q}\right]
.\\
\end{array}$$
Now consider  the rational numbers
$1-\frac{1}{p}+\frac{1}{q}$ and
$\frac{1}{p}-\frac{1}{q}$. Both numbers are $2$-adically
integral and belong to the open unit interval of the real line,
but  one and only one  of them (namely the former) is a
$2$-adic unit. Therefore one has
$\alpha_{pq}=-1$. 

\subsection{Calculation of $\alpha_{pq}$ in the case $2=p$}
Assume that $2=p$.
\subsubsection{}
One has 
$\cbold_{\sigma_2}=\cbold_1=0$ in this case and hence
$$\alpha_{2q}=\sign
\sigma_2\sin\cbold_{\sigma_q}.$$
 
\subsubsection{}
We may assume that the set $T$ has been chosen in such a way
that the lifting operator $H=\Hbold_T$ has the following
properties:
$$
\begin{array}{rcl}
H\sigma_q^i\left[\frac{4}{q}\right]&=&
\left\{\begin{array}{cl}
\sigma_q^i\left[\frac{4}{q}\right]\;\;\;\;\;\;\;&\mbox{for
$0\leq i<(q-1)/2$,}\\\\
 0\;\;\;\;\;\;\;&\mbox{for $(q-1)/2\leq i<q-1$,}
\end{array}\right.\\\\
H\left[\frac{\nu}{4}\right]&=&
\left\{\begin{array}{cl}
\left[\frac{\nu}{4}\right]\;\;\;\;\;\;\;\;\;\;\;&\mbox{for
$\nu=q$,}\\\\
0\;\;\;\;\;\;\;\;\;\;\;&\mbox{for
$\nu=-q$,}
\end{array}\right.\\\\
H\sigma_q^i\left[\frac{1}{2}+\frac{2}{q}\right]&=&
\left\{\begin{array}{cl}
\sigma_q^i\left[\frac{1}{2}+\frac{2}{q}\right]&
\mbox{for $0\leq i<(q-1)/2$,}\\\\
0&\mbox{for $(q-1)/2\leq i<q-1$,}
\end{array}\right.\\\\
H\sigma_q^i\left[\frac{\nu}{4}+\frac{1}{q}\right]&=&
\left\{\begin{array}{cl}
0&\mbox{for
$0\leq i<q-1$ and $\nu=1$,}\\\\
\sigma_q^i\left[\frac{\nu}{4}+\frac{1}{q}\right]&\mbox{for
$0\leq i<q-1$ and $\nu=-1$.}\end{array}\right.\\\\
\end{array}
$$

\subsubsection{}
Put
$$\Theta_q:=\left(\sigma_q^0+
\cdots+\sigma_q^{(q-3)/2}\right).$$
The relations
$$
(1-\sigma_q)\Theta_q[a+b]
=\left(1-\sigma_q^{\frac{q-1}{2}}\right)[a+b]=
[a+b]-[a-b]\;\;\;\;\;\;
\left(a\in \frac{1}{4}\ZZ,\;b\in
\frac{1}{q}\ZZ/\ZZ\right)$$ 
are crucial to our calculations.

\subsubsection{}
Recall that in the case $p=2$ one has
$$\xbold=[0]+\left[\frac{1}{2}\right].$$
We have
$$
\begin{array}{rcl}
Y_q\xbold&=&-\Theta_q
\left(\left[\frac{4}{q}\right]+\left[-\frac{4}{q}\right]
+\left[\frac{1}{2}+\frac{2}{q}\right]
+\left[\frac{1}{2}-\frac{2}{q}\right]\right)\\\\
HY_q\xbold&=&
-\Theta_q
\left(\left[\frac{4}{q}\right]+\left[\frac{1}{2}+\frac{2}{q}\right]
\right),\\\\
H(1-\sigma_q)
HY_q\xbold
&=&-\left(\left[\frac{4}{q}\right]+
\left[\frac{1}{2}+\frac{2}{q}\right]\right)\\\\
Y_2H(1-\sigma_q)
HY_q\xbold
&=&-\left(\left[\frac{4}{q}\right]-\left[\frac{2}{q}\right]
-\left[\frac{1}{4}+\frac{1}{q}\right]
-\left[-\frac{1}{4}+\frac{1}{q}\right]\right)\\\\
Y_2\xbold&=&-\left[\frac{q}{4}\right]-\left[-\frac{q}{4}\right]\\\\
HY_2\xbold&=&
-\left[\frac{q}{4}\right],\\\\
(Y_2HY_q-Y_q
HY_2)\xbold&=&
-\Theta_q
\left(\left[\frac{4}{q}\right]-\left[\frac{2}{q}\right]
-\left[\frac{1}{4}+\frac{1}{q}\right]
-\left[-\frac{1}{4}+\frac{1}{q}\right]\right)
\\\\
&&+\left[\frac{q}{4}\right]-\left[\frac{1}{4}\right]
-\Theta_q
\left(\left[\frac{1}{4}+\frac{1}{q}\right]
+\left[\frac{1}{4}-\frac{1}{q}\right]\right),\\\\
\abold&=&
\Theta_q
\left[-\frac{1}{4}+\frac{1}{q}\right]-H\Theta_q
\left(\left[\frac{4}{q}\right]-
\left[\frac{2}{q}\right]\right)+
H\left(\left[\frac{q}{4}\right]-
\left[\frac{1}{4}\right]\right),\\\\
H(1-\sigma_q)\abold&
=&
\left[-\frac{1}{4}+\frac{1}{q}\right]
-\left[-\frac{1}{4}-\frac{1}{q}\right]
-H(1-\sigma_q)H\Theta_q
\left(
\left[\frac{4}{q}\right]-
\left[\frac{2}{q}\right]\right),\\\\
HY_2H(1-\sigma_q)
HY_q\xbold
&=&\left[-\frac{1}{4}+\frac{1}{q}\right]-H\left(\left[\frac{4}{q}\right]-
\left[\frac{2}{q}\right]\right)
\\\\

\cbold_{\sigma_q}&=&H(1-\sigma_q)\abold-
HY_2H(1-\sigma_q)HY_q\xbold+HY_q
H(1-\sigma_q)HY_2\xbold
\\\\&=&-
\left[-\frac{1}{4}-\frac{1}{q}\right]+
H(1-(1-\sigma_q)H\Theta_q)
\left(\left[\frac{4}{q}\right]-
\left[\frac{2}{q}\right]\right).
\end{array}
$$
\subsubsection{}
Finally, one has
$$\frac{1}{\sin\cbold_{\sigma_q}}\equiv
\sin\pi\left(\frac{1}{4}+
\frac{1}{q}\right)\equiv
\frac{1}{\sqrt{2}}
\bmod{
\left(\RR\cap
\QQ\left(\ebold\left(\frac{1}{4q}\right)\right)
\right)^\times}
$$
because the sines and cosines of all integral multiples of
the angle $\pi/q$ belong to
the
field
$\RR\cap\QQ\left(\ebold\left(\frac{1}{4q}\right)\right)$,
and hence
$\alpha_{2q}=
\sign \sigma_2\sqrt{2}=-1$.
The proof of the Main Formula is
complete.

\section{Acknowledgements}
I thank P.~Das, Y.~Ouyang and S.~Seo for insightful
conversations and correspondence on topics related to
this paper. I also thank D.~Thakur for very helpful
criticism of a first draft of this paper.

\end{document}